\documentclass[reqno,12pt]{amsart}
\usepackage{mathrsfs}
\usepackage{amssymb,latexsym}
\usepackage{amsmath,amsfonts,amsthm}
\usepackage{graphicx} 
\usepackage{indentfirst}
\usepackage{fancyhdr}
\usepackage{CJK}
\usepackage[margin=2.5cm]{geometry}
\usepackage{color}
\usepackage{txfonts}

\def\doteqdot{:=}

\begin{document}

\newtheorem{theorem}{Theorem}[section]
\newtheorem{corollary}[theorem]{Corollary}  
\newtheorem{lemma}[theorem]{Lemma}
\newtheorem{definition}[theorem]{Definition}
\newtheorem{proposition}[theorem]{Proposition}
\newtheorem{remark}[theorem]{Remark}
\newtheorem{example}[theorem]{Example}
\newtheorem{theoremalph}{Theorem}
 \renewcommand\thetheoremalph{\Alph{theoremalph}}

\newcommand{\bta}{\begin{theoremalph}}
\newcommand{\eeta}{\end{theoremalph}}
\newcommand{\bth}{\begin{theorem}}
\newcommand{\eeth}{\end{theorem}}
\newcommand{\ble}{\begin{lemma}}
\newcommand{\ele}{\end{lemma}}
\newcommand{\bco}{\begin{corollary}}
\newcommand{\eco}{\end{corollary}}
\newcommand{\bde}{\begin{definition}}
\newcommand{\ede}{\end{definition}}
\newcommand{\bpr}{\begin{proposition}}
\newcommand{\epr}{\end{proposition}}
\newcommand{\bre}{\begin{remark}}
\newcommand{\ere}{\end{remark}}
\newcommand{\beg}{\begin{example}}
\newcommand{\eeg}{\end{example}}

\newcommand{\beq}{\begin{equation}}
\newcommand{\eeq}{\end{equation}}
\newcommand{\ben}{\begin{equation*}}
\newcommand{\een}{\end{equation*}}
\newcommand{\beqn}{\begin{eqnarray}}
\newcommand{\eeqn}{\end{eqnarray}}
\newcommand{\be}{\begin{eqnarray*}}
\newcommand{\ee}{\end{eqnarray*}}
\newcommand{\ban}{\begin{align*}}
\newcommand{\ean}{\end{align*}}
\newcommand{\bal}{\begin{align}}
\newcommand{\eal}{\end{align}}
\newcommand{\bad}{\aligned}
\newcommand{\ead}{\endaligned}
\newcommand{\lan}{\langle}
\newcommand{\ran}{\rangle}

\newcommand{\na}{\nabla}
\newcommand{\vp}{\varphi}
\newcommand{\La}{\Lambda}
\newcommand{\la}{\lambda}
\newcommand{\Om}{\Omega}
\newcommand{\ta}{\theta}
\newcommand{\fr}{\frac}
\newcommand{\iy}{\infty}
\newcommand{\ve}{\varepsilon}
\newcommand{\pa}{\partial}
\newcommand{\al}{\alpha}
\newcommand{\mr}{\mathbb{R}^n}
\newcommand{\bu}{\bullet}
\newcommand{\si}{\sigma}
\newcommand{\he}{{\rm Hess}}
\newcommand{\rc}{{\rm Ric}}
\newcommand{\dv}{{\rm div}}

\title
{Lower Bound Estimates for The First Eigenvalue of The Weighted $p$-Laplacian on Smooth Metric Measure Spaces}
\author{Yu-Zhao Wang}
\author{Huai-qian Li$^1$}
\footnote{Partially supported by the National Natural Science Foundation of China (NSFC) No.11401403 and the Australian Research Council (ARC) grant DP130101302. }
\footnote{\it The paper has been accepted by Differential Geometry and its Applications.}
\address{Institute of Applied Mathematics, Academy of Mathematics and Systems Science \\
Chinese Academy of Sciences, Beijing 100190, China}
\email{yuzhao@amss.ac.cn}
\address{School of Mathematics, Sichuan University, Chengdu 610064, China;
Department of Mathematics, Macquarie University, NSW 2109 Sydney, Australia}
\email{huaiqianlee@gmail.com; huaiqian.li@mq.edu.au}

\maketitle
\maketitle \numberwithin{equation}{section}
\maketitle \numberwithin{theorem}{section}
\setcounter{tocdepth}{2}
\setcounter{secnumdepth}{2}

\begin{abstract}
New lower bounds of the first nonzero eigenvalue of the weighted $p$-Laplacian are established on compact smooth metric measure spaces with or without boundaries. Under the assumption of positive lower bound for the $m$-Bakry--\'{E}mery Ricci curvature, the Escober--Lichnerowicz--Reilly type estimates are proved; under the assumption of nonnegative $\infty$-Bakry--\'{E}mery Ricci curvature and the $m$-Bakry--\'{E}mery Ricci curvature bounded from below by a non-positive constant, the Li--Yau type lower bound estimates are given. The weighted $p$-Bochner formula and the weighted $p$-Reilly formula are derived as the key tools for the establishment of the above results.

\vspace{5mm}
\vspace{5mm}

\end{abstract}
\section{Introduction and main results}
Let $(M,g)$ be a compact $n$-dimensional Riemannian manifold. Denote by ${\rm Ric}$ the Ricci curvature, by $D$ the diameter of $M$ and by $\partial M$ the boundary of $M$. It is known that the eigenvalue estimate for the Laplacian on Riemannian manifolds, which is studied intensively, is an important and long-standing issue  in geometric analysis and PDE theory (see e.g. \cite{Ch,SY,PLi}). Among the results on the lowest eigenvalue estimates, the most famous one is obtained by Lichnerowicz \cite{Lichnerowicz} and Obata \cite{Obata}. Let $\lambda_1$ be the first nonzero eigenvalue of $M$.
\bta[Lichnerowicz--Obata] Let $(M,g)$ be a closed Riemannian manifold with positive Ricci curvature, i.e., ${\rm Ric}\ge Kg$ for some $K>0$. Then
\begin{equation}\label{Lich}
\lambda_1\ge \frac{n}{n-1}K,
\end{equation}
and equality holds if and only if $M$ is isometric to the $n$-sphere with constant sectional curvature $\frac{K}{n-1}$.
\eeta

When $M$ has a nonempty boundary, if ${\rm Ric}\ge Kg$ for some $K>0$ and the mean curvature with respect to the outward unit norm vector field is nonnegative, Reilly \cite[Theorem 4]{Reilly} established the Lichnerowicz--Obata result for the first Dirichlet eigenvalue, and if the boundary of $M$ is convex, Escobar \cite[Theorem 4.3]{Escobar} proved the Lichnerowicz--Obata result for the first nonzero Neumann eigenvalue.

Assume that $M$ is closed and has nonnegative Ricci curvature, i.e., ${\rm Ric}\ge 0$. Li--Yau \cite{LiYau} deduced that
\begin{eqnarray}\label{LY}
\lambda_1\geq \frac{\pi^2}{2D^2}.
\end{eqnarray}
Later, Zhong--Yang \cite{ZhongYang} improved this estimate by an optimal lower bound, i.e.,
\begin{eqnarray}\label{ZY}
\lambda_1\ge \left(\frac{\pi}D\right)^2.
\end{eqnarray}
Moreover, Hang--Wang \cite{HW} proved that if the equality in \eqref{ZY} holds then  $M$ must be isometric to the unit circle with radius $D/\pi$.

The results for the case of general Ricci curvature lower bound $K\in\mathbb{R}$ were obtained by Kr\"{o}ger \cite{Kr} using the gradient comparison technique, and independently by Chen--Wang \cite{ChenWang94,ChenWang97} applying the probabilistic ``coupling method'', and recently  by Andrews and Clutterbuck \cite{AC2013} deriving sharp estimates on the  modulus of continuity for solutions of the heat equation. In a submitted paper, F.Z. Gong, D.J. Luo and the second named author \cite{GLL} also adapted successfully the coupling method to give a new proof of the fundamental gap conjecture, which was first proved by Andrews and Clutterbuck in the seminal paper \cite{AC2011}.

We should mention that there are also extensions of the Lichnerowicz--Obata theorem in the non-smooth metric measure spaces. Let $K>0$. Lott--Vinalli \cite{LV07} obtained \eqref{Lich} under the curvature-dimension condition ${\rm{CD}}(K,n)$ (see also \cite[Theorem 30.25]{Villani}). Erbar--Kuwada--Sturm \cite{EKS2013} showed \eqref{Lich} under the Reimannian curvature-dimension condition ${\rm{RCD}}^\ast(K,n)$, which is the strengthening of the curvature-dimension condition in the sense of Lott--Sturm--Villani by requiring the linearity of the heat flow, and Ketterer \cite{Ketterer14} established the rigidity part recently. In the Alexandrov space with a new notion of the Ricci curvature lower bound, which is stronger than the Riemannian curvature-dimension condition, namely $C(K,n)$, Qian-Zhang-Zhu \cite{QZZ2013} obtained the the Lichnerowicz--Obata theorem.

A natural question is how to estimate the bound of the first eigenvalue of the $p$-Laplacian on the compact Riemannian manifold $M$. Let $p\in (1,\infty)$ here and in the sequel. The $p$-Laplacian is defined by
$$\Delta_{p}u:={\rm div}(|\nabla u|^{p-2}\nabla u),\quad\mbox{for every }u\in W^{1,p}(M),$$
which is understood in distribution sense. The first nonzero eigenvalue of $\Delta_{p}$ is conventionally denoted by $\lambda_{1,p}$. When the closed Riemannian manifold $(M,g)$ has positive Ricci curvature lower bound, i.e., ${\rm Ric}\geq (n-1)Kg$ with $K>0$, by using a Bochner type formula,  Matei \cite{Mat} showed that
\begin{equation}\label{Matei}
\lambda_{1,p}\ge \left(\frac{(n-1)K}{p-1}\right)^{\frac{p}{2}},\quad p\geq 2,
\end{equation}
and also a Lichnerowicz--Obata type result by using the L\'{e}vy--Gromov isoperimetric inequality. On a compact Riemannian manifold with nonnegative Ricci curvature, Valtorta \cite[Theorem 7.1]{Val1} recently obtained a sharp lower bound for the first nonzero eigenvalue of  the $p$-Laplacian by applying the gradient comparison method, which can be viewed as a generalization of the Zhong--Yang estimate \eqref{ZY} when $p=2$, and he also established the rigidity part which we omit here (see \cite[Theorem 8.2]{Val1}).
\bta[Valtorta]\label{Val}
Let~$M$~be~{an} $n$-dimensional compact Riemannian manifold with nonnegative Ricci~curvature and possibly with convex boundary $\partial M$. Suppose that $\lambda_{1,p}$~{is} the first nonzero eigenvalue of the $p$-Laplacian (with the Neumann boundary condition if $\partial M\neq \emptyset$). Then
\beq\label{p-lower}
\lambda_{1,p}\ge (p-1)\left(\frac{\pi_p}D\right)^p,
\eeq
where\beq\label{pip}
\pi_p:=2\int^1_0\frac{ds}{(1-s^p)^{1/p}}=\frac{2\pi}{p\sin{(\pi/p)}}.
\eeq
\eeta

Now recall that a smooth metric measure space is a triple $(M,g, d\mu)$, where $(M,g)$ is a complete $n$-dimensional Riemannian manifold and $d\mu\doteqdot e^{-f}\,{dV}$ with $f$ a fixed smooth real-valued function on $M$. Denote by $\nabla$, $\Delta$ and $\rm{Hess}$ the gradient, Laplace and Hessian operators, and by $dV$ the Riemannian volume measure. The smooth metric measure space carries a natural analog of the Ricci curvature, the so-called $m$-Bakry--\'Emery Ricci curvature, which is defined as
$$
{\rm Ric}_f^m\doteqdot{\rm Ric}+{\rm Hess}f-\frac{\nabla f\otimes \nabla f}{m-n},\quad (n< m\leq \infty).
$$
In particular, when $m=\infty$, ${\rm Ric}^{\infty}_f\doteqdot{\rm Ric}_f\doteqdot{\rm Ric}+{\rm Hess}f$ is the classical Bakry--\'Emery Ricci curvature, which is introduced by Bakry--\'Emery \cite{BakryEmery} in the study of diffusion processes (see also \cite{BGL} for a comprehensive introduction), and then it has been extensively investigated in the theory of the Ricci flow (for example, the gradient Ricci soliton equation is precisely ${\rm Ric}_f =\lambda g$ for some constant $\lambda$) (see e.g. \cite{CLN}). The case where $m=n$ is only defined when $f$ is a constant function. There is also a natural analog of the Laplacian, that is, the weighted Laplacian (also called the $f$-Laplacian, drifting Laplacian or Witten Laplacian in the literature), denoted by $\Delta_f=\Delta-\nabla f\cdot\nabla$, which is a self-adjoint operator in $L^2(M, d\mu)$.

For the first nonzero eigenvalue of the weighted Laplacian on the compact Riemannian manifold without boundary or with convex boundary, Bakry--Qian \cite[Theorem 14]{BQ} obtained a unified lower bound, under the assumption of the $m$-Bakry--\'{E}mery Ricci curvature bounded from below.  Under the inspiration of recent studies on Ricci solitons and self-shrinkers, there are many results including the ones on gradient estimates and eigenvalue estimates of the weighted Laplacian in smooth metric measure spaces via the $m$-Bakry--\'{E}mery Ricci curvature; see \cite{LiXD, LiXD2, MW1} and references therein.  Recently, under the assumption of  the Bakry--\'{E}mery Ricci curvature bounded from below, Futaki--Li--Li \cite{FLL} proved a lower bound for the first nonzero eigenvalue of the weighted Laplacian on the compact Riemannian manifold without boundary, i.e.,
\begin{equation*}
\lambda_{1}\ge\sup_{s\in(0,1)}\left\{4s(1-s)\frac{\pi^2}{D^2}+sK\right\}.
\end{equation*}
Li--Wei \cite[Theorem 3]{LW} also showed a Lichnerowicz--Obtata type lower bound of the first nontrivial eigenvalue of the weighted Laplacian on a compact Riemannian manifold with boundary under the assumption of positive $m$-Bakry--\'{E}mery Ricci curvature lower bound.

In this note, we consider the first eigenvalue lower bound estimate of the weighted $p$-Laplacian, denoted by $\Delta_{p,f}$, on the compact smooth metric measure space $(M,g,d\mu)$. More precisely, for a function $u\in W^{1,p}(M)$, the weighted $p$-Laplacian is defined by
$$\Delta_{p,f}u\doteqdot e^{f}{\rm div}(e^{-f}|\nabla u|^{p-2}\nabla u),$$
which is also understood in distribution sense. We call that $\lambda$ is an eigenvalue of the weighted $p$-Laplacian $\Delta_{p,f}$  if there exists a nonzero function $u\in W^{1,p}(M)$  satisfying
\begin{equation}\label{pweigen}
\Delta_{p,f}u=-\lambda |u|^{p-2}u
\end{equation}
in distribution sense. We also use the notation $\lambda_{1,p}$ to denote the first nonzero eigenvalue of  $\Delta_{p,f}$ (since we will never consider $\Delta_p$ from now on).

 Inspired by works mentioned above, combining the weighted $p$-Bochner~formula (see Lemma \ref{pbochner} below) and the weighted $p$-Reilly formula (see Lemma \ref{pReillyformula} below) with the gradient estimate technique, we can obtain some lower bound estimates for the eigenvalue $\lambda_{1,p}$ in terms of the sign of the $m$-Bakry--\'{E}mery Ricci curvature.

Now we begin to introduce the main results in this work. We first need some notations. Let $\mathbf{n}$ be the outer unit normal vector field of $\partial M$. The second fundamental form of $\partial M$ is defined by ${\rm II}(X,Y)=\langle\nabla_X \mathbf{n}, Y\rangle$ for any vector fields $X$ and $Y$ on $\partial M$. The quantities
$$H(x)\doteqdot {\rm tr}({\rm II}_x)\quad\mbox{and}\quad H_f(x)\doteqdot H(x)-\langle\nabla f(x),\mathbf{n}(x)\rangle$$
are the mean curvature and the weighted mean curvature of $\partial M$ at $x\in M$. We call that $\partial M$ is convex if the second fundamental form ${\rm II}\ge0$. If $\partial M\neq \emptyset$, we assume the Dirichlet boundary condition or Neumann boundary condition, and then denote the first Dirichlet eigenvalue and the first nonzero Neumann eigenvalue of the weighted $p$-Laplacian  $\Delta_{p,f}$ as $\lambda_{Dir}$ and  $\lambda_{Neu}$, respectively.
\bth[Escobar--Lichnerowicz--Reilly type estimates]\label{LichEst} Assume $p\ge2$ and $K>0$. Let $(M,g,d\mu)$ be a compact smooth metric measure space.\\
{\rm (1)} Suppose $\partial M=\emptyset$. If ${\rm Ric}^m_f\ge Kg$, then,
\beq\label{Lichnerowicz}
\lambda_{1,p}\ge\frac{1}{(p-1)^{p-1}}\left(\frac {mK}{m-1}\right)^{\frac p2},
\eeq
where $\lambda_{1,p}$ is the first nonzero eigenvalue of the weighted $p$-Laplacian. In particular, if ${\rm Ric}_f\ge Kg$,
\beq\label{Lichnerowicz2}
\lambda_{1,p}\ge\frac{K^{\frac p2}}{(p-1)^{p-1}}.
\eeq
{\rm(2)} Suppose $\partial M\neq\emptyset$.  For the Dirichlet eigenvalue $\lambda_{Dir}$, we assume the weighted mean curvature $H_f$ is nonnegative, and for the Neumann eigenvalue $\lambda_{Neu}$, we assume the boundary $\partial M$ is convex. If ${\rm Ric}^m_f\ge Kg$, then
\beq\label{Lichnerowicz3}
\lambda_{Dir}\ge\frac{1}{(p-1)^{p-1}}\left(\frac {mK}{m-1}\right)^{\frac p2},\quad \lambda_{Neu}\ge\frac{1}{(p-1)^{p-1}}\left(\frac {mK}{m-1}\right)^{\frac p2}.
\eeq
In particular, if ${\rm Ric}_f\ge Kg$, then
\beq\label{LichReilly2}
\lambda_{Dir}\ge\frac{K^{\frac p2}}{(p-1)^{p-1}},\quad \lambda_{Neu}\ge\frac{K^{\frac p2}}{(p-1)^{p-1}}.
\eeq
\eeth
\bre
(i) It seems that even for the $p$-Laplacian case, estimates in \eqref{Lichnerowicz}-\eqref{LichReilly2} are new. In particular, when $p=2$, namely the weighted Laplacian case, the lower bounds in \eqref{Lichnerowicz} and \eqref{Lichnerowicz3} reduce to $\frac m{m-1}K$, which was obtained in \cite[Theorem 3]{LW} (see also \cite[Theorem 2]{MD}), and in addition, if $f$ is constant, then for $\partial M=\emptyset$, \eqref{Lichnerowicz} is exactly the Lichnerowicz estimate \eqref{Lich}, and for $\partial M\neq\emptyset$, estimates in \eqref{Lichnerowicz3} are due to \cite{Escobar} and independently \cite{Xia} in the Dirichlet case and to \cite{Reilly} in the Neumann case.

(ii) Compared with the lower bound in \cite[Theorem 3.2]{Mat} for the particular $p$-Laplacian case (see also \eqref{Matei}), it is easy to see that the lower bound in \eqref{Lichnerowicz2} is strictly smaller when $p>2$.

(iii) Recently, under the assumption of Theorem \ref{LichEst} (1), L.F. Wang \cite[Theorem 1.1]{LFW} obtained an estimate analogous to \eqref{Lichnerowicz} for the case where $m<\infty$; however, it is quite involved compared with our result.
\ere

The next two results generalize the ones in the Laplacian setting due to Li--Yau \cite{LiYau} (see also \cite[Chapter 3]{SY} and \eqref{LY} above).
\bth\label{peigest} Let $(M,g,d\mu)$ be a compact smooth metric measure space with nonnegative Barky--\'Emery Ricci curvature (i.e., ${\rm Ric}_f\geq 0$).\\
{\rm (1)} Suppose $\partial M=\emptyset$. Then the first nonzero eigenvalue of the weighted $p$-Laplacian satisfies
\begin{equation*}\label{peigen}
\lambda_{1,p}\ge (p-1)\left(\frac{\pi_p}{2D}\right)^p.
\end{equation*}\\
{\rm (2)} Suppose $\partial M\neq \emptyset$. If the weighted mean curvature $H_f$ of $\partial M$ is nonnegative, then
\begin{equation*}\label{Dpeigen}
\lambda_{Dir}\ge (p-1)\left(\frac{\pi_p}{2D}\right)^p;
\end{equation*}
if $\partial M$ is convex, then
\begin{equation*}\label{Npeigen}
\lambda_{Neu}\ge (p-1)\left(\frac{\pi_p}{2D}\right)^p.
\end{equation*}
Here $\pi_p$ is defined in \eqref{pip}.
\eeth
\bre
In \cite{Val1,NV}, the authors obtained sharp estimates for the first nonzero eigenvalue of the $p$-Laplacian on compact Riemannian manifolds by a different approach (see also Theorem \ref{Val}). However, we cannot prove the sharp form under the condition of Theorem \ref{peigest} by our approach for the moment. We should mention that also for the $p$-Laplacian, under the assumption of the quasi-positive Ricci curvature (i.e., ${\rm{Ric}}\geq 0$ with at least one point where $\rm{Ric}>0$), the same lower bounds for $\lambda_{1,p}$ and $\lambda_{Dir}$ were obtained in \cite{zhanghc07}.
\ere
\bth\label{LowerEst}  Let $(M,g,d\mu)$ be a closed smooth metric measure space satisfying ${\rm Ric}^m_f\ge-Kg$ with $K\geq 0$. Then the first nonzero eigenvalue of the weighted $p$-Laplacian satisfies
\begin{equation}\label{LiYau}
{\lambda_{1,p}}\ge\frac{C(p,m)}{D^p}\exp\Big(-\sqrt{(m-1)K}D\Big),
\end{equation}
where
$$C(p,m)=\frac2{m+1}\Big(\frac p{p-1}\Big)^{p-1}e^{-p}.$$
In addition, if $M$ has a convex boundary, then the first nonzero Neumann eigenvalue $\lambda_{Neu}$ also satisfies the estimate \eqref{LiYau}.
\eeth
\begin{remark}
In the setting of the $p$-Laplacian, a lower bound estimate similar to \eqref{LiYau} is established in \cite[Theorem 4]{LFW09} without mentioning the nonempty boundary case.
\end{remark}

This note is organized as follows. In Section 2, we establish the weighted $p$-Bochner and the weighted $p$-Reilly formulas as the key tools to prove our main results by a detailed calculation. The proofs of main theorems are presented in Section 3.

\section*{Acknowledgements}
The first named author would like to thank Professor Xiang-Dong Li for his interest and their illuminating discussion. Part of the work was done when the second named author was visiting Professor Thierry Coulhon at MSI in The Australian National University from 8 October 2014 to 5 February 2015, and he would like to thank his invitation. The two authors would also like to thank Professor Huichun Zhang for sending them the interesting paper \cite{zhanghc07}. And thanks are also given to an anonymous referee for his or her careful reading and comments, which make the presentation of this paper better.

\section{The weighted $p$-Bochner and  $p$-Reilly formulas}
Let $(M,g,d\mu)$ be a smooth metric measure space.  In the Laplacian setting, the important tools for gradient estimates and first eigenvalue estimates are the maximum principle, Bochner formula and Reilly formula. Likewise, in our weighted $p$-Laplacian setting, we are going to establish the weighted $p$-Bochner formula and weighted $p$-Reilly formula. Instead of using the Laplacian, we use the linearized operator of the weighted $p$-Laplacian at the maximum point in this work. The linearized operator of the weighted $p$-Laplacian at point $u\in C^2(M)$ such that $\nabla u\neq0$ is given by (see \cite{CW,WJC})
\ben\bad
\mathcal{L}_f(\psi)\doteqdot& e^f{\rm div}\left(e^{-f}|\nabla u|^{p-2}A(\nabla \psi)\right)\\
=&|\nabla u|^{p-2}\Delta_f\psi+(p-2)|\nabla u|^{p-4}{\rm Hess}\,\psi(\nabla u,\nabla u)\\
&+(p-2)\Delta_{p,f}u\frac{\langle\nabla u,\nabla\psi\rangle}{|\nabla u|^2}+2(p-2)|\nabla u|^{p-4}{\rm Hess}\, u\left(\nabla u,\nabla\psi-\frac{\nabla u}{|\nabla u|}\Big\langle\frac{\nabla u}{|\nabla u|},\nabla\psi\Big\rangle\right),
\ead\een
for a smooth function $\psi$ on $M$, where $A$ can be viewed as a tensor and it is defined as
\begin{equation}\label{Atensor}
A\doteqdot {\rm{Id}}+(p-2)\frac{\nabla u\otimes \nabla u}{|\nabla u|^2}.
\end{equation}
This operator $\mathcal{L}_f$ is defined pointwise only at the points that $\nabla u\neq0$ holds, and moreover, at these points, $\mathcal{L}_f$ is strictly elliptic. Denote by $\mathrm{L}_f$ the sum of the second order parts of $\mathcal{L}_f$, and hence
\ben
\mathrm{L}_f(\psi)=|\nabla u|^{p-2}\Delta_f\psi+(p-2)|\nabla u|^{p-4}{\rm Hess}\, \psi(\nabla u,\nabla u).
\een
Note that $$\mathcal{L}_f(u)=(p-1)\Delta_{p,f}u,\quad\mathrm{L}_f(u)=\Delta_{p,f}u.$$

We should mention that, since the equation \eqref{pweigen} can be either degenerate or singular  at the points such that $\nabla u=0$ (according to the value of $p$),  we usually use an $\varepsilon$-regularization technique by replacing the linearized operator $\mathcal{L}_f$ with its approximate operator, i.e.,
$$\mathcal{L}_{f,{\varepsilon}}\psi=e^f{\rm div}\left(e^{-f} w_{\varepsilon}^{\frac p2-1}A_{\varepsilon}(\nabla \psi)\right),$$
where $\varepsilon>0$, $w_{\varepsilon}=|\nabla u_{\varepsilon}|^2+\varepsilon$, $A_{\varepsilon}=\rm{Id}+(p-2)\frac{\nabla u_{\varepsilon}\otimes \nabla u_{\varepsilon}}{w_{\varepsilon}}$ and $u_{\varepsilon}$ is a solution to the approximate equation
\ben
e^f{\rm div}\big(e^{-f} w_{\varepsilon}^{\frac p2-1}\nabla u_{\varepsilon}\big)=-\lambda|u_{\varepsilon}|^{p-2}u_{\varepsilon}.
\een
In order to avoid the tedious presentation, we omit the details here; the interested reader should refer to \cite{KoNi} for example.

Now we  prove a nonlinear form of the weighted Bochner type formula related to the linearized weighted $p$-Laplacian.
\ble[weighted $p$-Bochner formula]\label{pbochner}
Let $(M,g,d\mu)$ be a smooth metric measure space. Given a $C^3$ function $u$, if $|\nabla u|\neq0$, then
\beq\bad\label{pBoch3}
\frac1p\mathrm{L}_f(|\nabla u|^p)=&|\nabla u|^{2p-4}\Big(|{\rm Hess}\,u|^2+p(p-2)(\Delta_{\infty}u)^2+{\rm Ric}_f(\nabla u,\nabla u)\Big)\\
&+|\nabla u|^{p-2}\Big(\langle\nabla\Delta_{p,f}u,\nabla u\rangle-(p-2)\Delta_{\infty}u\Delta_{p,f}u\Big),
\ead\eeq
where $\Delta_{\infty}u:=\frac{{\rm Hess}\,u(\nabla u,\nabla u)}{|\nabla u|^2}$;
furthermore,
\begin{align}\label{pBoch4}
\frac1p\mathcal{L}_f(|\nabla u|^p)
=|\nabla u|^{2p-4}\Big(|{\rm Hess}\,u|_A^2+{\rm Ric}_f(\nabla u,\nabla u)\Big)+|\nabla u|^{p-2}\langle\nabla u,\nabla\Delta_{p,f}u\rangle,
\end{align}
where $|{\rm Hess}\, u|^2_A=A^{ik}A^{jl}u_{ij}u_{kl}$ and $A$ is defined in \eqref{Atensor}.
\ele
\begin{proof}
Choose a local orthonormal frame of vector fields $e_1,\cdots,e_n$ with $e_n=\mathbf{n}$ on the boundary $\partial M$, and denote $u_i=du(e_i)$, $u_{ij}={\rm Hess}\,u(e_i,e_j)$, etc. Let $w\doteqdot|\nabla u|^2$. Then $w_i=2u_ku_{ki}$, $w_j=2u_lu_{lj}$ and $w_{ij}=2u_{kj}u_{ki}+2u_ku_{kij}$, so that $\Delta_{\infty}u=\frac{u_iu_{ij}u_j}w=\frac{\langle\nabla w,\nabla u\rangle}{2w}$.
Let us calculate $\Delta_f(|\nabla u|^p)$ and ${\rm Hess}(|\nabla u|^p)$. Applying the weighted Bochner formula (see e.g. \cite{LiXD} or \cite{WW}), we have
\ben\bad
\frac1p\Delta_f(w^{p/2})=&\frac12w^{p/2-1}\left(\Delta_fw+(\frac p2-1)\frac{|\nabla w|^2}w\right)\\
=&w^{p/2-1}\left(|{\rm Hess}\,u|^2+{\rm Ric}_f(\nabla u,\nabla u)+\langle\nabla u,\nabla\Delta_fu\rangle+(p-2)\frac{|\nabla w|^2}{4w}\right),
\ead\een
and
\begin{align*}
\frac1p\nabla_i\nabla_j(w^{p/2})=&\frac12w^{p/2-1}\left(w_{ij}+(\frac p2-1)\frac{w_iw_j}{w}\right)\\
=&w^{p/2-1}\left(u_{kj}u_{ki}+u_ku_{kij}+(p-2)\frac{u_ku_{ki}u_lu_{lj}}{w}\right),
\end{align*}
where $\nabla_i=\nabla_{e_i}$.  These lead us to
\ben\bad
\frac1p\frac{(\nabla_i\nabla_jw^{p/2})u_iu_j}{w}=&w^{p/2-1}\left(\frac{u_{kj}u_ju_{ki}u_i}w+\frac{u_ku_iu_ju_{kij}}w+(p-2)\frac{u_ku_{ki}u_iu_lu_{lj}u_j}{w^2}\right)\\
=&w^{p/2-1}\left(\frac{|\nabla w|^2}{4w}+\frac{u_ku_iu_ju_{kij}}w+(p-2)(\Delta_{\infty}u)^2\right).
\ead\een
Then we obtain
\ben\bad
\frac1p\mathrm{L}_f(|\nabla u|^p)=&\frac1pw^{p/2-1}\left(\Delta_fw^{p/2}+(p-2)\frac{(\nabla_i\nabla_jw^{p/2})u_iu_j}{w}\right)\\
=&w^{p-2}\left(|{\rm Hess}\,u|^2+{\rm Ric}_f(\nabla u,\nabla u)+\langle\nabla u,\nabla\Delta_fu\rangle+(p-2)\frac{|\nabla w|^2}{4w}\right)\\
&+(p-2)w^{p-2}\left(\frac{|\nabla w|^2}{4w}+\frac{u_ku_iu_ju_{kij}}w+(p-2)(\Delta_{\infty}u)^2\right).
\ead\een
On the other hand,
\begin{align*}
\langle\nabla u,\nabla\Delta_{\infty}u\rangle=&u_k\cdot\left(\frac{u_iu_ju_{ij}}{w}\right)_k\\
=&\frac{u_{ijk}u_iu_j u_k}{w}+\frac{u_ku_{ik}u_ju_{ij}}{w}+\frac{u_ku_{jk}u_iu_{ij}}{w}-\frac{2u_ku_lu_{kl}u_iu_ju_{ij}}{w^2}\\
=&\frac{u_{ijk}u_iu_j u_k}{w}+\frac{|\nabla w|^2}{2w}-2(\Delta_{\infty}u)^2,
\end{align*}
and
\begin{align*}
&\langle\nabla u,\nabla\Delta_{p,f}u\rangle=\left\langle\nabla u,\nabla\Big(w^{p/2-1}(\Delta_fu+(p-2)\Delta_{\infty}u)\Big)\right\rangle\\
=&w^{p/2-1}\Big(\langle\nabla u,\nabla\Delta_fu\rangle+(p-2)\langle\nabla u,\nabla\Delta_{\infty}u\rangle+(\frac p2-1)w^{-p/2}\langle\nabla u,\nabla w\rangle\Delta_{p,f}u\Big)\\
=&w^{p/2-1}\langle\nabla u,\nabla\Delta_fu\rangle
+(p-2)w^{p/2-1}\Big(\frac{u_{ijk}u_iu_j u_k}{w}+\frac{|\nabla w|^2}{2w}-2(\Delta_{\infty}u)^2+w^{1-p/2}\Delta_{\infty}u\Delta_{p,f}u\Big).
\end{align*}
Thus, we conclude that
\ben\bad
\frac1p\mathrm{L}_f(|\nabla u|^p)&=w^{p-2}\Big(|{\rm Hess}\,u|^2+{\rm Ric}_f(\nabla u,\nabla u)+w^{1-p/2}\langle\nabla u,\nabla\Delta_{p,f}u\rangle\Big)\\
&+(p-2)w^{p-2}\Delta_{\infty}u\Big(p\Delta_{\infty}u-w^{1-p/2}\Delta_{p,f}u\Big).
\ead\een
Thus, \eqref{pBoch3} follows.

Moreover,  the first order part of $\mathcal{L}_f$ is given by
\begin{align*}
&\frac1p\mathcal{L}_f(w^{p/2})-\frac1p{\mathrm{L}_f}(w^{p/2})\\
=&(p-2)w^{p/2-1}\Delta_{p,f}u\frac{\langle\nabla u,\nabla w\rangle}{2w}+(p-2)w^{p-3}{\rm Hess}\, u\left(\nabla u,\nabla w-\frac{\langle\nabla u,\nabla w\rangle}{w}\nabla u\right)\\
=&(p-2)w^{p/2-1}\Delta_{p,f}u\Delta_{\infty}u+(p-2)w^{p-2}\left(\frac{|\nabla w|^2}{2w}-2(\Delta_{\infty}u)^2\right).
\end{align*}
Thus,
\begin{align*}
\frac1p\mathcal{L}_f(|\nabla u|^p)=&w^{p-2}\Big(|{\rm Hess}\,u|^2+{\rm Ric}_f(\nabla u,\nabla u)+w^{1-p/2}\langle\nabla u,\nabla\Delta_{p,f}u\rangle\Big)\\
&+(p-2)^2w^{p-2}(\Delta_{\infty}u)^2+(p-2)w^{p-2}\frac{|\nabla w|^2}{2w}\\
=&w^{p-2}\Big({\rm Ric}_f(\nabla u,\nabla u)+w^{1-p/2}\langle\nabla u,\nabla\Delta_{p,f}u\rangle\Big)\\
&+w^{p-2}\Big(|{\rm Hess}\, u|^2+\frac{p-2}2\frac{|\nabla w|^2}w+\frac{(p-2)^2}4\frac{\langle\nabla u,\nabla w\rangle^2}{w^2}\Big)\\
=&w^{p-2}\Big(|{\rm Hess}\,u|_A^2+{\rm Ric}_f(\nabla u,\nabla u)+w^{1-p/2}\langle\nabla u,\nabla\Delta_{p,f}u\rangle\Big),
\end{align*}
where we used the fact that
$|{\rm Hess}\,u|_A^2=|{\rm Hess}\, u|^2+\frac{p-2}2\frac{|\nabla w|^2}w+\frac{(p-2)^2}4\frac{\langle\nabla u,\nabla w\rangle^2}{w^2}$. Therefore, the proof of \eqref{pBoch4} is complete.
\end{proof}

To understand the case when $M$ has a nonempty boundary, we need an identity which follows from integration by parts on the weighted $p$-Bochner formula with respect to $d\mu$. We call it a weighted $p$-Reilly formula and present it in the next theorem. Denote by $d\sigma=e^{-f}dV_{\partial M}$ the weighted Riemannian volume measure of $\partial M$, where $dV_{\partial M}$ is the Riemannian volume measure of $\partial M$. Let $\nabla_{\partial}$ and $\Delta_{\partial}$ be the covariant derivative and Laplacian on $\partial M$  with respect to the induced Riemannian metric, and let $\Delta_{\partial,f}$ be the weighted Laplacian on $\partial M$.

\begin{theorem}[weighted $p$-Reilly formula]\label{pReillyformula}
Let $(M,g,d\mu)$ be a compact smooth metric measure space with boundary $\partial M$. Then for any $C^3$ function $u$,
\begin{align}\label{PReilly}
&\int_M\left[(\Delta_{p,f}u)^2-|\nabla u|^{2p-4}\Big(|{\rm Hess}\, u|_A^2+{\rm Ric}_f(\nabla u,\nabla u)\Big)\right]\,d\mu\notag\\
=&\int_{\partial M}|\nabla u|^{2p-4}\Big[(H_fu_n+\Delta_{\partial,f}u)u_n+{\rm II}(\nabla_{\partial} u,\nabla_{\partial} u)-\langle\nabla_{\partial} u,\nabla_{\partial} u_{n}\rangle\Big]\,d\sigma,
\end{align}
where $u_n=\frac{\partial u}{\partial\mathbf{n}}$.
\end{theorem}
\begin{proof} We use the same notations in the proof of Lemma \ref{pbochner}. Integrating the $p$-Bochner formula \eqref{pBoch4} over the manifold $M$ with respect to  $d\mu=e^{-f}\,{dV}$, we have
\ben
\int_M\Big(\frac1p\mathcal{L}_f(|\nabla u|^p)-|\nabla u|^{p-2}\langle\nabla u,\nabla(\Delta_{p,f}u\rangle\Big)\, d\mu
=\int_M|\nabla u|^{2p-4}\Big(|{\rm Hess}\, u|_A^2+{\rm Ric}_f(\nabla u,\nabla u)\Big)\,d\mu.
\een
Integration by parts immediately yields
\begin{align*}
&\int_M|\nabla u|^{p-2}\langle\nabla u,\nabla\Delta_{p,f}u\rangle\, d\mu
=\int_{\partial M}|\nabla u|^{p-2} (\Delta_{p,f}u)u_n\,d\sigma-\int_M(\Delta_{p,f}u)^2\,d\mu\\
=&\int_{\partial M}|\nabla u|^{2p-4} \sum_{i,j=1}^n\left(u_{ii}+(p-2)\frac{u_iu_ju_{ij}}{|\nabla u|^2}-u_if_i\right)u_n\,d\sigma-\int_M(\Delta_{p,f}u)^2\,d\mu.
\end{align*}
On the other hand, the divergence theorem implies
\ben
\frac1p\int_M\mathcal{L}_f(|\nabla u|^p)\,d\mu=
\int_{\partial M}|\nabla u|^{2p-4}\sum_{i,j=1}^n\left(u_iu_{in}+(p-2)\frac{u_iu_ju_{ij}u_n}{|\nabla u|^2}\right)\,d\sigma.
\een
Thus,
\begin{align*}
&\int_M\Big(\frac1p\mathcal{L}_f(|\nabla u|^p)-|\nabla u|^{p-2}\langle\nabla u,\nabla\Delta_{p,f}u\rangle\Big)\, d\mu\\
=&\int_{\partial M}|\nabla u|^{2p-4}\sum_{i=1}^n\left(u_iu_{in}-u_{ii}u_n+u_if_iu_n\right)\,d\sigma+\int_M(\Delta_{p,f}u)^2\,d\mu\\
=&\int_{\partial M}|\nabla u|^{2p-4}\left(\sum_{i=1}^{n-1}(u_iu_{in}-u_{ii}u_n+u_if_iu_n)+f_nu^2_n\right)\,d\sigma+\int_M(\Delta_{p,f}u)^2\,d\mu
\end{align*}
Following the similar calculation in \cite{chwang} or \cite{MD}, we have
\begin{align*}
\sum_{i=1}^{n-1}u_{ii}
=&\sum_{i=1}^{n-1}\big(e_i(e_iu)-(\nabla_{e_i}e_i) u\big)\\
=&\sum_{i=1}^{n-1}\big({\nabla_\partial}_{e_i}e_i-\nabla_{e_i}e_i\big)u+\Delta_{\partial}u\\
=& Hu_n+\Delta_{\partial}u,
\end{align*}
and
\begin{align*}
\sum_{i=1}^{n-1}u_{in}u_i=&\sum_{i=1}^{n-1}u_{ni}u_i =\sum_{i=1}^{n-1}(e_i(e_nu)-(\nabla_{e_i}e_n)u)u_i\\
=&\sum_{i=1}^{n-1}e_i(u_n)u_i-\sum_{i,j=1}^{n-1}\rm{II}_{ij}u_iu_j\\
=&\langle \nabla_\partial u,\nabla_\partial u_n\rangle-\rm{II}(\nabla_\partial u,\nabla_\partial u).
\end{align*}

Combining all these identities, we obtain
\ben\bad
&\int_M\left[(\Delta_{p,f}u)^2-|\nabla u|^{2p-4}\Big(|{\rm Hess}\, u|_A^2+{\rm Ric}_f(\nabla u,\nabla u)\Big)\right]\,d\mu\\
=&\int_{\partial M}|\nabla u|^{2p-4}\left[\Big(H_fu_n+\Delta_{\partial,f}u\Big)u_n+\Big({\rm II}(\nabla_{\partial} u,\nabla_{\partial} u)-\langle\nabla_{\partial} u,\nabla_{\partial} u_n\rangle\Big)\right]\,d\sigma.
\ead\een
This finishes the proof of \eqref{PReilly}.
\end{proof}
\bre
{\rm (1)} For the general Riemannian manifold $(M,g)$, the $p$-Reilly formula \eqref{PReilly} is also new. In fact, when $f=const.$, the identity \eqref{PReilly} reduces to
\begin{align*}\label{PReilly2}
&\int_M\left[(\Delta_{p}u)^2-|\nabla u|^{2p-4}\Big(|{\rm Hess}\, u|_A^2+{\rm Ric}(\nabla u,\nabla u)\Big)\right]\,dV\\
=&\int_{\partial M}|\nabla u|^{2p-4}\Big[(Hu_n+\Delta_{\partial}u)u_n+{\rm II}(\nabla_{\partial} u,\nabla_{\partial} u)-\langle\nabla_{\partial} u,\nabla_{\partial} u_{n}\rangle\Big]\,dV_{\partial M}.
\end{align*}
{\rm(2)} For $p=2$, \eqref{PReilly} becomes the classic Reilly formula for the weighted Laplacian $\Delta_f$ (see \cite{MD,KM}), i.e.,
\begin{equation*}\label{PReilly3}
\int_M\left[(\Delta_{f}u)^2-\Big(|{\rm Hess} u|^2+{\rm Ric}_f(\nabla u,\nabla u)\Big)\right]\,d\mu
=\int_{\partial M}\Big[(H_fu_n+2\Delta_{\partial,f}u)u_n+{\rm II}(\nabla_{\partial} u,\nabla_{\partial} u)\rangle\Big]\,d\sigma.
\end{equation*}
\ere

\section{The first eigenvalue estimate for weighted $p$-Laplacian}

In this section, we prove the main results, namely Theorems \ref{LichEst}, \ref{peigest} and \ref{LowerEst}, which are presented in the following three subsections, respectively.

\subsection{Positive $m$-Bakry--\'{E}mery Ricci curvature}

The original idea of the proof of Theorem \ref{LichEst} comes from the Laplacian case (see e.g. Appendix A in \cite{CLN}). The main tools we use here are the weighted $p$-Bochner formula \eqref{pBoch4} and the weighted $p$-Reilly formula \eqref{PReilly}, and some tricks from the geometry related to the Bakry--\'{E}mery Ricci curvature (see \cite{BGL,LiXD,WW}).

\begin{proof}[Proof of Thoerem \ref{LichEst}]
Let $K>0$ and let $(M,g,d\mu)$ be a compact smooth metric measure space with ${\rm Ric}^m_f\ge Kg$. Since the case where $p=2$ is obviously true, we should assume that $p>2$.

For $\lambda\in \mathbb{R}$, let $u$ be the solution to the  equation \eqref{pweigen} and let $m<\infty$. The weighted $p$-Bochner formula \eqref{pBoch4} implies that
\begin{align}\label{peigenest1}
\frac1p\mathcal{L}_f(|\nabla u|^p)=&|\nabla u|^{2p-4}(|{\rm Hess}\, u|_A^2+{\rm Ric}_f(\nabla u,\nabla u))+\langle|\nabla u|^{p-2}\nabla u,\nabla\Delta_{p,f}u\rangle\notag\\
\ge&\frac{(\Delta_{p,f}u)^2}{m}+|\nabla u|^{2p-4}{\rm Ric}^m_f(\nabla u,\nabla u)+\langle|\nabla u|^{p-2}\nabla u,\nabla\Delta_{p,f}u\rangle\notag\\
\ge&\frac{\lambda^2}{m}|u|^{2p-2}+K|\nabla u|^{2p-2}-(p-1)\lambda|\nabla u|^{p}|u|^{p-2},
\end{align}
where the first inequality follows from (see e.g. \cite{LiXD})
\begin{align}\label{Cauchy}
&|\nabla u|^{2p-4}\Big(|{\rm Hess}\, u|_A^2+{\rm Ric}_f(\nabla u,\nabla u)\Big)\notag\\
\ge&\frac1n\left(|\nabla u|^{p-2}{\rm tr}_A({\rm Hess}\, u)\right)^2+|\nabla u|^{2p-4}{\rm Ric}_f(\nabla u,\nabla u)\notag\\
=&\frac1n(\Delta_{p,f}u+|\nabla u|^{p-2}\langle\nabla u,\nabla f\rangle)^2+|\nabla u|^{2p-4}\left({\rm Ric}^m_f(\nabla u,\nabla u)+\frac{\langle\nabla u,\nabla f\rangle^2}{m-n}\right)\notag\\
\ge&\frac1m(\Delta_{p,f}u)^2+|\nabla u|^{2p-4}{\rm Ric}^m_f(\nabla u,\nabla u),
\end{align}
and the second inequality follows from \eqref{pweigen} and the assumption.\\
(1). Suppose $\partial M=\emptyset$. Integrating \eqref{peigenest1} on $M$ with respect to measure $d\mu=e^{-f}\,{dV}$, we have
\begin{equation}\label{Lich1}
0\ge\frac{\lambda^2}{m}\int_M|u|^{2p-2}\,d\mu+K\int_M|\nabla u|^{2p-2}\,d\mu-(p-1)\lambda\int_M|\nabla u|^{p}|u|^{p-2}\,d\mu.
\end{equation}
Multiplying  $\Delta_{p,f}u+\lambda|u|^{p-2}u=0$ by $|u|^{p-2}u$ on both sides and integrating by parts, we obtain
\beq\label{Lich2}
\lambda\int_M|u|^{2p-2}\,d\mu=(p-1)\int_M|\nabla u|^p|u|^{p-2}\,d\mu.
\eeq
Using the H\"older inequality, for any $p>2$, we have
\begin{align*}
\lambda\int_M|u|^{2p-2}\,d\mu
\le&(p-1)\left(\int_M\big(|\nabla u|^p\big)^{\frac{2p-2}{p}}\,d\mu\right)^{\frac{p}{2p-2}}
\left(\int_M\big(|u|^{p-2}\big)^{\frac{2p-2}{p-2}}\,d\mu\right)^{\frac{p-2}{2p-2}}\\
=&(p-1)\left(\int_M|\nabla u|^{2p-2}\,d\mu\right)^{\frac{p}{2p-2}}
\left(\int_M|u|^{2p-2}\,d\mu\right)^{\frac{p-2}{2p-2}};
\end{align*}
hence
\begin{align}\label{Lich4}
\int_M|\nabla u|^{2p-2}\,d\mu\ge\left(\frac{\lambda}{p-1}\right)^{2-\frac2p}\int_M|u|^{2p-2}\,d\mu.
\end{align}
Combining \eqref{Lich1}, \eqref{Lich2} and \eqref{Lich4}, we obtain
\ben
\left[\Big(\frac1m-1\Big)\lambda^2+K\left(\frac{\lambda}{p-1}\right)^{2-\frac2p}\right]\int_M|u|^{2p-2}\,d\mu\le0.
\een
Since $\int_M|u|^{2p-2}\,d\mu>0$, for any $p>2$,  we have
\ben
\lambda\ge\frac{1}{(p-1)^{p-1}}\left(\frac{mK}{m-1}\right)^{\frac p2},
\een
which is  \eqref{Lichnerowicz} when $m<\infty$.

For the case $m=\infty$, if ${\rm Ric}_f\ge K$, then by \eqref{pBoch4} and \eqref{Lich2}, we get
\ben
\left[-\lambda^2+K\left(\frac{\lambda}{p-1}\right)^{2-\frac2p}\right]\int_M|u|^{2p-2}\,d\mu\le0,
\een
and hence \eqref{Lichnerowicz2} follows.\\
(2). Suppose  $(M,g,d\mu)$ has a nonempty boundary $\partial M$. Using the weighted $p$-Reilly  formula \eqref{PReilly} and the estimate in \eqref{Cauchy}, we have
\begin{align*}
&\int_M\left[\left(1-\frac1m\right)\lambda^2|u|^{2p-2}-K|\nabla u|^{2p-2}\right]\,d\mu\notag\\
\ge&\int_{\partial M}|\nabla u|^{2p-4}\Big[H_fu^2_n+{\rm II}(\nabla_{\partial} u,\nabla_{\partial} u)+(\Delta_{\partial,f}u)u_n-\langle\nabla_{\partial} u,\nabla_{\partial} u_n\rangle\Big]\,d\sigma.
\end{align*}
Applying the technique used for estimating \eqref{Lich2} and \eqref{Lich4}, we get
\begin{align*}\label{pDN}
&-\left[\Big(\frac1m-1\Big)\lambda^2+K\left(\frac{\lambda}{p-1}\right)^{2-\frac2p}\right]\int_M|u|^{2p-2}\,d\mu\notag\\
\ge&\int_{\partial M}|\nabla u|^{2p-4}\Big[H_fu^2_n+{\rm II}(\nabla_{\partial} u,\nabla_{\partial} u)+(\Delta_{\partial,f}u)u_n-\langle\nabla_{\partial} u,\nabla_{\partial} u_n\rangle\Big]\,d\sigma.
\end{align*}
Recall the Dirichlet boundary condition is $u=0$ on $\partial M$ and the Neumann boundary condition is $\frac{\partial u}{\partial \mathbf{n}}=0$ on $\partial M$.  Using the assumption on the boundary $\partial M$, we finally obtain the two estimates \eqref{Lichnerowicz3} for $m<\infty$. In a similar spirit, we also have \eqref{LichReilly2}.
\end{proof}

\subsection{Nonnegative Bakry--\'{E}mery curvature}

In order to prove  Theorem \ref{peigest}, we compute the linearized weighted $p$-Laplacian $\mathrm{L}_f$  of a proper gradient quantity, and then apply the maximum principle and the weighted $p$-Bochner formula \eqref{pBoch3}.

Let $(M,g,d\mu)$ be a compact smooth metric measure space without boundary. Assume $u$ is the eigenfunction corresponding to $\lambda_{1,p}$, which means $u$ satisfies the equation \eqref{pweigen} with $\lambda$ replaced by $\lambda_{1,p}$ in distribution sense. It is known that the solution $u$ is  $C^{1,\alpha}(M)$ for some $0<\alpha<1$. Then $\int_M|u|^{p-2}u\, d\mu=0$ implies that $u$ changes the sign. Thus, multiplying it by some constant, we can assume that
\ben\label{eigen1}
\max_{x\in M} u(x)=1-\varepsilon\quad\mbox{and}\quad\min_{x\in M}u(x)=-a(1-\varepsilon),
\een
where $\varepsilon$ is a small enough positive number and $0 < a\le 1$. Then $u$ satisfies the following gradient estimate.
\begin{lemma}\label{gradient}
Let $(M,g,d\mu)$ be a closed smooth metric measure space and let $u$ be the eigenfunction corresponding to the eigenvalue $\lambda_{1,p}$ satisfying \eqref{eigen1}. Assume the Bakry-\'{E}mery Ricci curvature is nonnegative, i.e., ${\rm Ric}_f\geq 0$. Then it holds
\begin{equation}\label{eigengrad}
\frac{|\nabla u|^p}{1-|u|^p}\le\frac{\lambda_{1,p}}{p-1}.
\end{equation}
\end{lemma}
\begin{proof}
Let $F\doteqdot\frac{|\nabla u|^p}{1-|u|^p}$. Suppose it achieves its maximum at a point $x_0 \in M$. If $|\nabla u(x_0)|=0$, then \eqref{eigengrad} trivially holds.  If $|\nabla u(x_0)|\neq0$, then we may rotate the frame so that $u_1(x_0)=|\nabla u(x_0)|$ and $u_j(x_0)=0$ for $j\ge2$ and $F$ is $C^\infty$ in a neighborhood of $x_0$. By the maximum principle, we have
\ben
\nabla F(x_0)=0,\quad\mathrm{L}_fF(x_0)\le0.
\een
At the point $x_0$, we obtain
\ben\label{eigen2}
u_{11}=-u|u|^{p-2}u_1^{2-p}F,\quad u_{1j}=0,\;\;j\ge2.
\een
Let $\lambda=\lambda_{1,p}$. Combining this with \eqref{pweigen} and the weighted $p$-Bochner formula \eqref{pBoch3}, we deduce that
\beq\label{eigen3}
\frac1p\mathrm{L}_f(|u|^p)=(p-1)^2|u|^{p-2}|\nabla u|^p-\lambda |u|^{2p-2},
\eeq
and
\begin{eqnarray}\label{eigen4}
\frac1p\mathrm{L}_f(|\nabla u|^p)&=&|\nabla u|^{2p-4}\Big(|{\rm Hess}\,u|^2+{\rm Ric}_f(\nabla u,\nabla u)+p(p-2)(\Delta_{\infty}u)^2\Big)\notag\\
&&+\lambda|\nabla u|^{p-2}|u|^{p-2}\Big((p-2)u\Delta_{\infty}u-(p-1)|\nabla u|^2\Big)\\
&\ge&(p-1)^2u_1^{2p-4}u^2_{11}+\lambda u_1^{p-2}|u|^{p-2}\Big((p-2)u_{11}u - (p-1)u_1^2\Big)\notag,
\end{eqnarray}
where we have used ${\rm Ric}_f\ge0$, $|{\rm Hess}\, u|^2\ge u^2_{11}$ and $\Delta_{\infty}u\ge u_{11}$.

Applying $\mathrm{L}_f$ to both sides of the identity $(1-|u|^p)F=|\nabla u|^p$, with the fact $\nabla F(x_0)=0$ and $\mathrm{L}_fF(x_0)\le0$, we get
\ben
\mathrm{L}_f(|\nabla u|^p)\le-F\mathrm{L}_f(|u|^p).
\een

Combining this with \eqref{eigen3} and \eqref{eigen4}, we obtain
\begin{align*}
&(p-1)^2u_1^{2p-4}u^2_{11}+\lambda u_1^{p-2}|u|^{p-2}\Big((p-2)u_{11}u-(p-1)u_1^2\Big)\\
\le&- F \Big((p-1)^2|u|^{p-2}u_1^p-\lambda |u|^{2p-2}\Big).
\end{align*}
Noticing that $u_1^{p-2}u_{11}=-u|u|^{p-2}F$ and $(1-|u|^p)F=|\nabla u|^p$, we have
\begin{align*}
&(p-1)^2|u|^{2p-2}F^2-\lambda |u|^{p-2}\Big((p-2)|u|^pF+(p-1)F(1-|u|^p)\Big)\\
\le&-F\Big((p-1)^2|u|^{p-2}F(1-|u|^p)-\lambda |u|^{2p-2}\Big),
\end{align*}
which immediately implies
\begin{equation*}
F\le\frac{\lambda}{p-1}.
\end{equation*}
Therefore, we finish the proof.
\end{proof}

\begin{proof}[Proof of Theorem \ref{peigest}]
(1) Let $x_1, x_2 \in M$ such that $u(x_1) = 1-\varepsilon$ and $u(x_2) = -a(1-\varepsilon)$. Take a normalized minimal geodesic from $x_2$ to $x_1$. Since the Riemannian distance between $x_1$ and $x_2$ is no bigger than $D$, \eqref{eigengrad} implies that
\ben
D\Big(\frac{\lambda}{p-1}\Big)^{\frac1p} \ge\int^{1-\varepsilon}_{-a(1-\varepsilon)}\frac{du}{(1-|u|^p)^{1/p}}.
\een
Letting $\varepsilon\to0^+$, we get
\ben
D\Big(\frac{\lambda}{p-1}\Big)^{\frac1p}\ge\int^{1}_{-a}\frac{du}{(1-|u|^p)^{1/p}} > \int^{1}_{0}\frac{du}{(1-u^p)^{1/p}}
=\frac{\pi}{p\sin{(\pi/p)}}=\frac{\pi_p}2,
\een
which concludes the proof of Theorem \ref{peigest} (1).

(2)  Assume that $\partial M$ is nonempty. Let $F:=\frac{|\nabla u|^p}{1-|u|^p}$ be the auxiliary function as in Lemma \ref{gradient}. If the maximum of $F$ is attained in the interior of $M$, then using the similar argument as in (1), we can still obtain the desired estimate. Thus, we assume the maximum of $F$ is attained at, say $x_0\in \partial M$, and then derive a contradiction.

Choosing  an orthonormal frame of vector fields $\{e_i\}^n_{i=1}$ in a neighborhood of $x_0$ in $\partial M$ such that $e_n=\mathbf{n}$. By the Hopf maximum principle, we have
\ben
F_n(x_0)\doteqdot\frac{\partial F}{\partial \mathbf{n}}(x_0)>0.
\een
Taking the derivative $\frac{\partial }{\partial \mathbf{n}}$ to both sides of the identity $|\nabla u|^p=(1-|u|^p)F$, we get
\begin{align*}
p|\nabla u|^{p-2}\sum^n_{i=1}u_iu_{in}=F_n(1-|u|^p)-pF|u|^{p-2}uu_n.
\end{align*}

In either the case of the Dirichlet problem ($u=0$ on $\partial M$) or the Neumann problem ($u_n=0$), we have
\begin{align}\label{DirNeu}
\frac1p|\nabla u|^{2-p}(1-|u|^p)F_n=\sum^n_{i=1}u_iu_{in}
=u_nu_{nn}+\sum^{n-1}_{i=1}u_iu_{in}.
\end{align}
For the Dirichlet problem, $u=u_i=0$ on $\partial M$ for $i\neq n$. Then \eqref{pweigen} implies
\begin{align*}
0=-\lambda|u|^{p-2}u=\Delta_{p,f}u=&|\nabla u|^{p-2}\Big(\Delta_fu+(p-2)\frac{{
\rm Hess}\,u(\nabla u,\nabla u)}{|\nabla u|^2}\Big)\\
=&|\nabla u|^{p-2}(\Delta u-f_nu_n+(p-2)u_{nn})\\
=&|\nabla u|^{p-2}\Big((p-1)u_{nn}-f_nu_n+\sum^{n-1}_{i=1}u_{ii}\Big).
\end{align*}
Thus,
\begin{align*}
(p-1)u_{nn}=&f_nu_n-\sum^{n-1}_{i=1}u_{ii}
=f_nu_n-\sum^{n-1}_{i=1}\Big(e_i(e_iu)-(\nabla_{e_i}e_i)u\Big)\\
=&f_nu_n-\sum^{n-1}_{i=1}\langle\nabla_{e_i}e_i,\mathbf{n}\rangle u_n
=f_nu_n-Hu_n=-H_fu_n.
\end{align*}
For the Neumann problem, $u_n=0$. Thus in either of the cases,
\beq\label{Dirichlet1}
u_nu_{nn}=-\frac1{p-1}H_fu^2_n.
\eeq

For $i\neq n$,
\begin{align}\label{Neumann2}
\sum^{n-1}_{i=1}u_iu_{in}=\sum^{n-1}_{i=1}u_i\Big(e_i(u_n)-\sum^{n-1}_{j=1}{\rm II}_{ij}u_j\Big)=-\sum^{n-1}_{i,j=1}{\rm II}_{ij}u_j,
\end{align}
since either $u_i=0$ on $\partial M$ (Dirichlet) or $e_i(u_n)=0$ (Neumann).

Putting \eqref{Dirichlet1} and  \eqref{Neumann2} into \eqref{DirNeu}, we obtain
\beq
\frac1p|\nabla u|^{2-p}(1-|u|^p)F_n=-\frac1{p-1}H_fu^2_n-\sum^{n-1}_{i,j=1}{\rm II}_{ij}u_iu_j.
\eeq

Therefore, in either the Dirichlet problem ($u=u_i=0$ on $\partial M$) with $H_f\ge0$ or the Neumann problem ($u_n=0$) with ${\rm II}_{ij}\ge0$, {$i,j=1,\cdots,n-1$}, we have
\ben
\frac1p|\nabla u|^{2-p}(1-|u|^p)F_n\le0,
\een
which is impossible since $|u|<1$ and $F_n>0$. This finishes the proof of Theorem \ref{peigest} (2).
\end{proof}

\subsection{Non-positive $m$-Bakry--\'{E}mery Ricci curvature}

In order to prove Theorem \ref{LowerEst}, we derive  a lemma at first.
\ble\label{GE}
Let $(M,g,d\mu)$ be a compact smooth metric measure space with the $m$-Bakry--\'{E}mery Ricci curvature ${\rm Ric}^m_f\ge-Kg$
for some constant $K\ge0$. Let the function $u: M\rightarrow [-1,1]$ satisfy equation \eqref{pweigen}.
For some constant $a>1$, define
\ben
G\doteqdot(p-1)^p|\nabla\log(a+u)|^p.
\een
Then,
\begin{align}
&\mathcal{L}_f(G)-\frac{\alpha}{p}\frac{|\nabla G|^2}{G^{2/p}}-\left(\frac{2(p-1)}{m-1}\Big(1+\lambda u w^{-\frac p2}h\Big)-p\right)\langle\nabla v,\nabla G\rangle w^{\frac p2-1}\notag\\
\ge&\frac p{m-1}G^2-pKG^{2(p-1)/p}+\Big(\frac{2u}{m-1}-a\Big)p\lambda hG+\frac{p}{m-1}(\lambda uh)^2,
\end{align}
where $\alpha\doteqdot\min\{2(p-1),\frac{m(p-1)^2}{m-1}\}$, $v:=(p-1)\log(a+u)$, $w\doteqdot(p-1)^2|\nabla \log(a+u)|^2$ and $h\doteqdot(p-1)^{p-1}\frac{|u|^{p-2}}{(a+u)^{p-1}}$.
\ele
\proof
Set $v=(p-1)\log(a+u)$.  From equation \eqref{pweigen}, it is easy to see that $v$ satisfies
\beq\label{eigenv}
\Delta_{p,f}v=-|\nabla v|^p-(p-1)^{p-1}\frac{\lambda|u|^{p-2}u}{(a+u)^{p-1}}.
\eeq
Set $w=|\nabla v|^2$. Then in terms of $w$, \eqref{eigenv} has the equivalent form
\beq\label{WPLap3}
(\frac p2-1)w^{p/2-2}\langle\nabla w,\nabla v\rangle+w^{p/2-1}\Delta_f v=-w^{p/2}-(p-1)^{p-1}\frac{\lambda|u|^{p-2}u}{(a+u)^{p-1}}.
\eeq
By using  the weighted $p$-Bochner formula \eqref{pBoch4}, we have
\begin{align}\label{pBochv}
\mathcal{L}_f(|\nabla v|^p)
=p|\nabla v|^{2p-4}\Big(|{\rm Hess}\,v|_A^2+{\rm Ric}_f(\nabla v,\nabla v)+|\nabla v|^{2-p}\langle\nabla v,\nabla\Delta_{p,f}v\rangle\Big),
\end{align}
where $\mathcal{L}_f$ is the linearized operator of $\Delta_{p,f}$ at $v$. Substituting \eqref{eigenv} into \eqref{pBochv}, we get
\beq\label{PbochG}
\mathcal{L}_f(G)=pw^{p-2}\Big(|{\rm Hess}\,v|_A^2+{\rm Ric}_f(\nabla v,\nabla v)\Big)
-pw^{\frac p2-1}\langle\nabla v,\nabla G\rangle-pa\lambda hG.
\eeq
where $h=(p-1)^{p-1}\frac{|u|^{p-2}}{(a+u)^{p-1}}$.

Now we estimate  $|{\rm Hess}\, v|^2_A$. We only need to estimate it over the point where $w>0$. Choose  a local orthonormal frame $\{e_i\}_{i=1}^n$ near any such a given point so that $\nabla v=|\nabla v|e_1$. Then $w=v_1^2$, $w_1=2v_{i1}v_i=2v_{11}v_1$, and for $j\ge2$, $w_j=2v_{j1}v_1$. Hence, $2v_{j1}=\frac{w_j}{w^{1/2}}$. From \eqref{WPLap3}, we immediately deduce that
\begin{align}\label{eigenv2}
\sum_{j=2}^nv_{jj}&=-w-(\frac p2-1)\frac{w_1v_1}{w}+f_1v_1-v_{11}-\lambda u w^{1-\frac p2}h\notag\\
&=-w-(p-1)v_{11}+f_1v_1-\lambda u w^{1-\frac p2}h.
\end{align}
By the definition of $A$ and the Cauchy-Schwarz inequality, we have
\begin{align*}
|{\rm Hess}\, v|_A^2&=|{\rm Hess}\, v|^2+\frac{(p-2)^2}{4w^2}\langle\nabla v,\nabla w\rangle^2+\frac{p-2}{2w}|\nabla w|^2\\
&=\sum^n_{i,k=1}v^2_{ik}+(p-2)^2v^2_{11}+2(p-2)\sum^n_{k=1}v^2_{1k}\\
&=(p-1)^2v^2_{11}+2(p-1)\sum^n_{k=2}v^2_{1k}+\sum^n_{i,k=2}v^2_{ik}\\
&\ge (p-1)^2v^2_{11}+2(p-1)\sum^n_{k=2}v^2_{1k}+\frac1{n-1}\Big(\sum_{j=2}^nv_{jj}\Big)^2.
\end{align*}
Substituting \eqref{eigenv2} into the above inequality, we obtain
\begin{align}\label{eigenv3}
|{\rm Hess}\, v|_A^2
\ge&(p-1)^2v^2_{11}+2(p-1)\sum^n_{k=2}v^2_{1k}+\frac1{n-1}\left(w+(p-1)v_{11}-f_1v_1+\lambda u w^{1-\frac p2}h\right)^2\notag\\
\ge&(p-1)^2v^2_{11}+2(p-1)\sum^n_{k=2}v^2_{1k}+\frac1{m-1}\left(w+(p-1)v_{11}+\lambda u w^{1-\frac p2}h\right)^2-\frac{(f_1v_1)^2}{m-n}\notag\\
\ge&\alpha\sum^n_{k=1}v^2_{1k}+\frac1{m-1}\left(w+\lambda u w^{1-\frac p2}h\right)^2
+\frac{2(p-1)v_{11}}{m-1}\left(w+\lambda u w^{1-\frac p2}h\right)-\frac{(f_1v_1)^2}{m-n},
\end{align}
where $\alpha=\min\{2(p-1),\frac{m(p-1)^2}{m-1}\}$, and we applied the inequality $(a-b)^2\ge\frac{a^2}{1+\delta}-\frac{b^2}{\delta}$ with $\delta=\frac{m-n}{n-1}>0$. Substituting the identities
\begin{align*}
2wv_{11}=\langle\nabla v,\nabla w\rangle\quad\mbox{and}\quad
\sum^n_{j=1}v^2_{1j}=\frac14\frac{|\nabla w|^2}{w}
\end{align*}
into \eqref{eigenv3}, we obtain
\begin{eqnarray}\label{eigenv4}
|{\rm Hess}\, v|_A^2
&\ge&\frac{\alpha}4\frac{|\nabla w|^2}w+\frac1{m-1}\left(w+\lambda u w^{1-\frac p2}h\right)^2\cr
&&+\frac{p-1}{m-1}\langle\nabla v,\nabla w\rangle\left(1+\lambda u w^{-\frac p2}h\right)-\frac{(f_1v_1)^2}{m-n}.
\end{eqnarray}

Combining \eqref{PbochG} and \eqref{eigenv4}, by the assumption that ${\rm Ric}^m_f\ge-Kg$,  we have
\begin{eqnarray}\label{WPLap11}
\mathcal{L}_f(G)
&\ge& pw^{p-2}\left(\frac{\alpha}4\frac{|\nabla w|^2}w+\frac1{m-1}\left(w+\lambda u w^{1-\frac p2}h\right)^2+\frac{p-1}{m-1}{\langle\nabla v,\nabla w\rangle}\Big(1+\lambda u w^{-\frac p2}h\Big) \right)\notag\\
&&-pw^{p-2} {\rm Ric}_f^m(\nabla v,\nabla v) - pw^{\frac p2-1}\langle\nabla v,\nabla G\rangle-a p\lambda hG\notag\\
&\ge&\frac p{m-1}G^2+\frac{\alpha}{p}\frac{|\nabla G|^2}{G^{2/p}}-pKG^{2(p-1)/p}+\Big(\frac{2u}{m-1}-a\Big)p\lambda hG+\frac{p}{m-1}(\lambda uh)^2\notag\\
&&+\left(\frac{2(p-1)}{m-1}\Big(1+\lambda u w^{-\frac p2}h\Big)-p\right)\langle\nabla v,\nabla G\rangle w^{\frac p2-1}\notag,
\end{eqnarray}
which is the desired result.
\endproof
\begin{proof}[Proof of Theorem \ref{LowerEst}]
Let $u$ be a nonconstant eigenfunction satisfying \eqref{pweigen}. Using the fact
\ben
-\lambda\int_{M}|u|^{p-2}u\,d\mu=\int_M\Delta_{p,f}u\,d\mu,
\een
we know that $u$ must change the sign. Without loss of generality, set $\min_M u=-1$ and $\max_M  u\le1$.

Consider the function $$v=(p-1)\log(a+u)$$ for some constant $a>1$. Assume that $x_0\in M$ is a point such that the function $G=|\nabla v|^p$ reaches the maximum. By using the maximum principle and Lemma \ref{GE}, we obtain
\ben\bad
0\ge&\frac p{m-1}G^2-pKG^{2(p-1)/p}+\Big(\frac{2u}{m-1}-a\Big)p\lambda h G+\frac{p}{m-1}(\lambda uh)^2\\
\ge&\frac p{m-1}G^2-pKG^{2(p-1)/p}+\Big(\frac{2u}{m-1}-a\Big)p\lambda h G.
\ead\een
Combining this and Young's inequality, we get
\begin{align*}
0\ge& G-(m-1)KG^{1-\frac2p}+(2u-(m-1)a)\lambda h\\
\ge&G-{\frac{p-2}p}G-\frac 2p((m-1)K)^{\frac p2}+(2u-(m-1)a)\lambda h\\
=&\frac 2pG-\frac 2p((m-1)K)^{\frac p2}+(2u-(m-1)a)\lambda h,
\end{align*}
which implies that, for all $x \in M$,
\ben\bad
G(x)\le& G(x_0)\le((m-1)K)^{\frac p2}-\frac p2(2u(x_0)-(m-1)a)\lambda (p-1)^{p-1}\frac{|u(x_0)|^{p-2}}{(a+u(x_0))^{p-1}}\\
\le&((m-1)K)^{\frac p2}+\frac {p(p-1)^{p-1}}2\Big(\frac{a}{a-1}\Big)^{p-1}(m+1)\lambda.
\ead\een

Integrating $G^{\frac1p}=(p-1)|\nabla\log(a+u)|$ along a minimal geodesic $\gamma:\,[0,1]\to M$ joining  points $x_1$ and $x_2$ such that $u(x_1) =-1$ and $u(x_2)=\max_M u$, we obtain
\ben\bad
(p-1)\log\left(\frac{a}{a-1}\right)\le&(p-1)\log\left(\frac{a+\max u}{a-1}\right)\\
\le&\int^1_0(p-1)|\nabla\log(a+u(\gamma(s)))||\dot{\gamma}(s)|\,ds\\
\le& D\left(((m-1)K)^{\frac p2}+\frac {p(p-1)^{p-1}}2\Big(\frac{a}{a-1}\Big)^{p-1}(m+1)\lambda\right)^\frac1p\\
\le&D\sqrt{(m-1)K}+D(p-1)^{\frac{p-1}p}\Big(\frac{a}{a-1}\Big)^{\frac{p-1}p}(m+1)^\frac1p\Big(\frac {p}2\lambda\Big)^\frac1p,
\ead\een
for all $a>1$, where we used the inequality $(x+y)^{\frac1p}\le x^{\frac1p}+y^{\frac1p}$ for $x,y\geq0$ and $p\ge1$.

Setting $t=\log\Big(\frac{a}{a-1}\Big)$ with $a>1$, we have
\beq\label{eigenv5}
\Big(\frac {p}2\lambda\Big)^\frac1p\ge\frac1D\Big(\frac{p-1}{m+1}\Big)^{
\frac1p}\Big(t-\frac{1}{p-1}\sqrt {(m-1) K}D\Big)e^{-\frac{p-1}{p}t}.
\eeq
 Choosing  $$t=\frac p{p-1}+\frac{1}{p-1}\sqrt {(m-1) K}D$$ such that, as a function of $t$, the right hand side of \eqref{eigenv5} achieves the maximum, we finally obtain the estimate
$$
\lambda\ge\frac2{m+1}\Big(\frac p{p-1}\Big)^{p-1}\frac1{D^p}e^{-\big(p+\sqrt{(m-1) K}D\big)}.
$$
The proof of \eqref{LiYau} is finished by letting
$$C(p,m)=\frac2{m+1}\Big(\frac p{p-1}\Big)^{p-1}e^{-p}.$$

In addition, if $\partial M$ is nonempty and convex, then following similarly the proof of \cite[Corollary 5.8]{PLi} in the Laplacian case, we complete the proof.
\end{proof}

\end{document}